\documentclass[12pt,a4paper]{amsart}
\usepackage[margin=25mm]{geometry}

 \usepackage{amsmath,amsthm}
\usepackage{amssymb}
\usepackage{epsfig}
\usepackage[all]{xypic}
\usepackage[colorlinks]{hyperref}
\usepackage[numbers]{natbib}
\usepackage{graphicx}

\newtheorem{theorem}{Theorem}[section]

\theoremstyle{definition}
\newtheorem{definition}[theorem]{Definition}

\newcommand{\R}{\mathbb{R}}

\newcommand{\p}{\mathbb{P}}
\newcommand{\Span}[1]{\mathrm{Span}\left\langle#1\right\rangle}
\DeclareMathOperator{\LL}{LGr}

\DeclareMathOperator{\Sp}{\mathsf{Sp}}

\newcommand{\E}{\mathcal{E}}

\newcommand{\Th}{^\textrm{th}}
\newcommand{\St}{^\textrm{st}}

\newcommand{\Nd}{^\textrm{nd}}

\begin{document}

\subjclass{30,80}


\title[Completely exceptional PDEs]{An   introduction to completely exceptional  $2\Nd$ order scalar PDEs}

\author[G. Moreno]{Giovanni Moreno}
\address{Institute of Mathematics, Polish Academy of Sciences\\
\'Sniadeckich 8, 00-956 Warszawa, Poland\\
E-mail: gmoreno@impan.pl}
\maketitle
 
 \begin{abstract}
 In his 1954 paper about the initial value problem for 2D hyperbolic nonlinear PDEs, P. Lax declared  that he had ``a strong reason to believe'' that there must exist a well--defined class of ``not genuinely nonlinear'' nonlinear PDEs. In 1978 G. Boillat coined the term ``completely exceptional'' to denote it. In the case of $2\Nd$ order (nonlinear) PDEs, he also proved that this class reduces to  the class of  Monge--Amp\`ere equations. We review here, against a unified  geometric background, the notion of complete exceptionality, the definition of a Monge--Amp\`ere equation, and the interesting link between them.
\end{abstract}

\tableofcontents

\section*{Introduction}
 
 A function $F=F(u_{11},u_{12},u_{22})$ in the three variables $u_{11},u_{12},u_{22}$ is linear if and only if $F$ is a solution to the system
 \begin{equation}\label{eqEqBanLin}
 \frac{\partial^2 F}{\partial u_{ij}\partial u_{kl}}=0\, ,\quad \forall i,j,k,l=1,2\, ,
 \end{equation}
 of $2\Nd$ order PDEs. Such a remark immediately becomes less silly when one begins thinking of the variables $u_{11},u_{12},u_{22}$ as $2\Nd$ order \emph{formal  derivatives} of a function $u$ in two extra variables $x_1,x_2$. Such a   perspective allows us to reinterpret $F$   as the left--hand side of a (scalar) $2\Nd$ order PDE in $u=u(x_1,x_2)$.\par
 Accordingly, \eqref{eqEqBanLin} must be thought of as an ``equation imposed on equations''. That is, the totality of the solutions of \eqref{eqEqBanLin} represents (the left--hand sides of) the equations which constitute  a special class of $2\Nd$ order scalar PDEs---the linear ones.   
 By making \eqref{eqEqBanLin} totally symmetric in the indices $i,j,k,l$, we obtain a weaker condition, viz.
 \begin{equation}\label{eqDerSymbPrimordiale}
  \frac{\partial^2 F}{\partial u_{(ij}\partial u_{kl)}}=0\, ,\quad \forall i,j,k,l=1,2\, .
 \end{equation}
 Obviously, \eqref{eqDerSymbPrimordiale} is satisfied by all (the left--hand sides of) linear $2\Nd$ order scalar PDEs. Less obviously, yet still straightforwardly, there is more than just linear PDEs in the solution set of  \eqref{eqDerSymbPrimordiale}. The curious reader may check this on the manifestly nonlinear  $F=u_{11}u_{22}-u_{12}^2$.\par 
Usually  a  $2\Nd$ order scalar PDE $F=0$ is accompanied by some initial data. If one is interested in the so--obtained initial value problem, then the first question to answer is whether or not the initial data are characteristic. Here comes to help a key geometric gadget associated with $F$, namely   the \emph{(principal) symbol}
\begin{equation}\label{eqSIMBOLOprimordiale}
 S(F):= \frac{\partial  F}{\partial u_{ij} }\partial_i\odot\partial_j\, 
 \end{equation} 
 of $F$. 
 Loosely speaking, \eqref{eqSIMBOLOprimordiale} is   a symmetric tensor on a 2--dimensional space.\footnote{We prefer to leave this space unspecified   in order to contain the size of this introduction.}  By abusing the notation, we will call $S(F)$   a ``metric'' even though $S(F)$ is controvariant and it is  degenerate in the majority of the really interesting cases.  By line--hyperplane duality,\footnote{Throughout this paper we use $^\circ$ to denote the annihilator of a linear subspace and, in particular, to express line--hyperplane duality.} the null directions of this ``metric'', which are easily   computed, can be interpreted as tangent hyperplanes in the space of the independent variables $x_i$'s. These are the \emph{characteristic hyperplane} associated with   the (generally nonlinear) equation  $F=0$. In the  geometric framework for   (nonlinear) PDEs described below, characteristic hyperplanes\footnote{The literature on this subject is boundless. The chief reference is Chapter V of the book \cite{MR1083148} by Bryant et al., but this might prove hard to the novice. The paper \cite{MR2985508} perhaps provides a slenderer introduction to the subject. The reader may also have a look at \cite{MorenoCauchy,Vitagliano2013,2017arXiv170104930S}.} are precisely the tangent spaces to the characteristic initial data, i.e., those for which uniqueness of the solution to the initial value problem is not guaranteed.\par
The symbol $S(F)$ is an intrinsic feature of the PDE $F=0$, in the sense that the \emph{tensors} $S(F_1)$ and $S(F_2)$ are the same (up to a projective factor) as long as $F_1=0$ and $F_2=0$ are the same equation. But the \emph{components}  $\frac{\partial  F}{\partial u_{ij} }$ of $S(F)$ may change. Similarly, \eqref{eqDerSymbPrimordiale} can be thought of as the components of a rank--4 symmetric tensor
\begin{equation}\label{eqEQS2Fprimordiale}
S^2(F)=\frac{\partial^2 F}{\partial u_{(ij}\partial u_{kl)}}\partial_i\odot\partial_j\odot\partial_h\odot\partial_k \, .
\end{equation}
Let $F_1=0$ be a PDE. A surprising discovery, originally obtained by  G. Boillat \cite{MR1139843} in 1991  and further clarified   in 2017 \cite{MR3603758}, is that the \emph{tensorial} equation\footnote{Here  ``$\approx$'' means ``proportional to''.}
\begin{equation}\label{eqCEprimordiale}
S^2(F)\approx S(F)\, , 
\end{equation}
where $F$ is such that $F=0$ and $F_1=0$ are the same equation, 
 is \emph{intrinsically associated} to the PDE $F_1=0$ (and \emph{not} to   its particular left--hand side $F_1$).  The equation \eqref{eqCEprimordiale} is clearly satisfied by (the left--hand sides of) linear equations (indeed $S(F)=0$, and the factor of proportionality can be set to zero), but its entire set of solutions is much more large. \par
In case the reader is wondering why we introduced the condition $S^2(F)\approx S(F)$ by starting from the trivial equation \eqref{eqEqBanLin}, that is why we introduced a certain class of PDEs by weakening the condition of linearity, the answer is simple. Because the class of linear PDEs is not invariant under contactomorphism, and then a larger class must exist. In fact, the class defined by Boillat by imposing \eqref{eqCEprimordiale} is even larger that the closure of the class of the linear PDEs under the action of the group of all contactomorphism. The role of  conctactomorphism in this context is clarified below.\par
Even though P. Lax dealt with systems of quasi--linear $1\St$ order PDEs and he never used the term ``completely exceptional'' (introduced---to the author best knowledge---by G. Boillat and T. Ruggeri in 1978 \cite{MR0481580}), the class of PDEs satisfying \eqref{eqCEprimordiale}  was called ``completely exceptional'' \emph{in the sense of Lax} by Boillat himself, referring to Lax's 1954 paper \cite{MR0068093}. And to this terminology we shall stick. The condition $S^2(F)\approx S(F)$  will always be ``the condition of complete exceptionality'' for the PDE $F=0$.\par
%
Below we provide a  solid geometric background to the condition of complete exceptionality $S^2(F)\approx S(F)$. We also explain how the fact that $F$ satisfies $S^2(F)\approx S(F)$ reflects on the behaviour of the solutions to an  initial value problem associated with the PDE $F=0$. In this short note, the reader  will find an  answer to the below questions.
\begin{enumerate}
\item In which sense \eqref{eqCEprimordiale} is intrinsically associated to $F=0$? (see Section \ref{sec1}.)
\item What is the set of all the solutions of \eqref{eqCEprimordiale}? (see Section \ref{sec2}.)
\item If a solution of \eqref{eqCEprimordiale} is interpreted in its turn as a $2\Nd$ order PDE, then what makes the latter \emph{exceptional} amongst all $2\Nd$ order PDEs? (see Section \ref{sec3}.)
\end{enumerate}
These answers already exist scattered throughout the literature \cite{MR3603758, MR1139843,MR0481580,MR1292999,MR1461074,MR571041,MR1461074}. The purpose of this note is precisely that of arranging   them in a unified self--consistent and minimalistic way.

\section{Geometry of $2\Nd$ order PDEs and their characteristics}\label{sec1}

Partial differential equations, together with their solutions and initial data, can be conveniently formalised in terms of smooth manifolds and differential forms on them. This is the core of E. Cartan's pioneering work, which evolved into the modern theory of Exterior Differential Systems (EDS) \cite{MR1083148}.
\subsection{$1\St$ order (scalar, nonlinear) PDEs} 
In the present context we are interested in PDEs in $n$ independent and 1 dependent variable. These, taken together, can be understood as the local coordinates    $(x_1,\ldots,x_n,u)$ of an $(n+1)$--dimensional manifold. For the purpose of studying  $1\St$ order PDEs, we need to add more variables, say $u_1,\ldots,u_n$.  We thus obtain a $(2n+1)$--dimensional manifold, henceforth denoted by $M$, with local coordinates  $(x_1,\ldots,x_n,u,u_1,\ldots,u_n)$.\par 
So far $u$ is just another name for  a coordinate, and as such it does not carry any dependence upon $x_1,\ldots,x_n$ whatsoever. Therefore, even if they seem so, the coordinates $u_1,\ldots,u_n$ are \emph{not} the derivatives of $u$. And the hypersurface
\begin{equation}\label{eq1stOrdDE}
\E=\{F(x_1,\ldots,x_n,u,u_1,\ldots,u_n)=0\}\subset M
\end{equation}
cut out by a function $F\in C^\infty(M)$ is \emph{not} a $1\St$ order PDE.\par 
In order to recover the intuition which is still obviously missing in the picture, it is enough to introduce the \emph{contact form}
\begin{equation}\label{eqContForm}
\theta=du-\sum_{i=1}^nu_idx_i
\end{equation}
on $M$. The \emph{contact manifold} $(M,\theta)$ is an example of an EDS. An EDS is just a manifold equipped with a set of differential forms (only $\theta$ in our example). (Such a fundamental role   of contact manifolds in the geometric framework for   PDEs explains why a  class of PDEs  needs to be closed under the group of contactomorphisms.)\par 
The main concern in the theory of EDS is to study the so--called \emph{variety of integral elements}. An integral element is simply a tangent plane to $M$, such that all the forms which define the EDS, together with    their differentials, vanish on it. Usually one groups integral elements according to their dimension. For example, let us study the $n$--dimensional integral elements of $M$, i.e., the $n$--dimensional tangent planes to $M$ such that both $\theta$ and $d\theta$ vanish on them.\par 
Let $L$ be such an $n$--dimensional integral element, and let $p\in M$ be the point of $M$ the $n$--plane $L$ is tangent to. In order to have $\theta$ to vanish on $L$, we need to pick the generators of $L$ from the hyperplane $\ker\theta_p\subset T_pM$. Plainly,
\begin{equation}
\ker\theta_p=\Span{\partial_1|_p,\ldots,\partial_n|_p,\partial_{u_1}|_p,\ldots,\partial_{u_n}|_p}\, ,
\end{equation}
where the vector fields $\partial_i:=\partial_{x_i}+u_i\partial_u$ are often called \emph{total derivatives}. Now we can ask when an $n$--plane of the form
\begin{equation}\label{eqLagrPlane}
L=\Span{\partial_{i}|_p+\sum_{j=1}^nu_{ij}\partial_{u_i}|_p}
\end{equation}
(on which by construction $\theta$ vanishes) makes also $d\theta$ vanish. Observe that the new symbols $u_{ij}$ appearing in \eqref{eqLagrPlane} are just numeric coefficients. By taking the differential of \eqref{eqContForm} and by imposing $d\theta|_L\equiv 0$, we find the simple condition
\begin{equation}
u_{[ij]}=0\, .
\end{equation}
That is, the space of $n$--dimensional integral elements tangent to $p\in M$ is parametrised by $n\times n$ symmetric matrices. Hence, the totality of all $n$--dimensional integral elements (tangent to arbitrary points of $M$) form a set, henceforth denoted by  $M^{(1)}$, naturally fibered over $M$, with abstract  fibre $S^2\R^n$.\par 
Now we can finally recover the familiar understanding of a partial differential equation. An \emph{integral submanifold} of an EDS is a submanifold $U\subset M$ all whose tangent spaces are   integral elements. In our example, an $n$--dimensional submanifold $U$ is integral if and only if there exists a \emph{function} $f$, such that
\begin{equation}\label{eqUUnoEffe}
U=U^1_f:=\left\{\left(x_1,\ldots,x_n,u=f(x_1,\ldots,x_n),\ldots, u_i=\frac{\partial f}{\partial x_i}(x_1,\ldots,x_n),\ldots\right)\right\}\, .
\end{equation}
Then the hypersurface \eqref{eq1stOrdDE} can be correctly interpreted as a $1\St$ order PDE, in the sense that $f$ is a solution to $\E$ if and only if $U_f^1$ is \emph{contained} into $\E$. But all of this is just a paraphrase of the Darboux theorem on the structure of Legendrian submanifolds of a contact manifold \cite{MR2352610}. In the literature, the Legendrian submanifolds $U_f^1$ are often called   \emph{graphs of $1\St$ jets of functions} \cite{MR1670044}.

\subsection{$2\Nd$ order (scalar, nonlinear) PDEs} \label{subPDEs}

It is somewhat useful to refer to the integral submanifolds $U$ of the contact EDS $(M,\theta)$ as \emph{candidate solutions}. Indeed, thanks to \eqref{eqUUnoEffe},   candidate solutions are in (a local) one--to--one correspondence with   functions in $n$ variables. A candidate solution may be thought of as a solution of the trivial equation $0=0$; it becomes a solution of the (nontrivial) equation \eqref{eq1stOrdDE} only if it is contained into $\E$. In a sense, the whole machinery so far introduced  just allowed to rephrase in terms of a set--theoretical inclusion the property for a function and its $1\St$ derivatives to satisfy a certain  relation.\par 
Given a candidate solution $U\subset M$, observe that all its tangent $n$--planes are, by defintion, integral elements, that is, points of $M^{(1)}$. In other words, $U$ may be as well considered as (an $n$--dimensional) submanifold of $M^{(1)}$, if we identify each point $p\in U$ with the corresponding  tangent space $T_pU$. Usually this is formalised by introducing a new set
\begin{equation}
U^{(1)}:=\{ T_pU\mid p\in U\}\, ,
\end{equation}
manifestly identical to $U$, but contained this time into $M^{(1)}$. It is customary to denote $(U_f^{1})^{(1)}$ simply by $U_f^2$ and call it the \emph{graph of the $2\Nd$ jet of $f$}.\par 
A hypersurface $\E\subset M^{(1)}$ is called a \emph{(scalar, nonlinear) $2\Nd$ order PDE (in $n$ independent variables)}. Indeed, in view of \eqref{eqLagrPlane}, such an hypersurface can be (locally) represented\footnote{From now on the symbol $\E$ refers to a hypersurface in $ M^{(1)}$ as in \eqref{eqEQ2ndOrdPDE} and not to a hypersurface in $M$ as in  \eqref{eq1stOrdDE}.} as
\begin{equation}\label{eqEQ2ndOrdPDE}
\E=\{ F=F(x_1,\ldots,x_n,u,u_1,\ldots,u_n,\ldots,u_{ij},\ldots)=0\}\subset M^{(1)}\, .
\end{equation} 
Then, a candidate solution $U_f^2$ is contained into $\E$ if and only if the function $f$, together with its $1\St$ and $2\Nd$ derivatives, fulfils the relation given  by $F=0$. Once again, the familiar intuition of a ($2\Nd$ order PDE) has been recovered.\par 
So far, we have just recast the well--known notions of  PDEs and their solutions in terms of hypersurfaces in $M^{(1)}$ and Legendrian submanifolds of $M$, respectively. In order to see the first nontrivial implication of such a reinterpretation, we need to inspect the vertical geometry of the bundle $M^{(1)}\longrightarrow M$.

\subsection{Geometry of the Lagrangian Grassmannian $\LL_{n,2n}$}\label{subsecLG}
The reader may have noticed that formula \eqref{eqLagrPlane} is not entirely accurate, in the sense that not all the $n$--planes in $\ker\theta_p$ are of that form. However, \eqref{eqLagrPlane} was useful to find a local description of the fibre $M_p^{(1)}$, that is the set of all integral $n$--planes at $p\in M$. In fact, the whole of $M_p^{(1)}$ is a topologically nontrivial compactification of the linear space of $S^2\R^n$, known as the \emph{Lagrangian Grassmannian} and usually denoted by $\LL_{n,2n}$.\par 
Let us regard a symmetric matrix $A\in S^2\R^n$ as a linear map from $\R^{n\ast}$ to $\R^n$, and let us extend it to a linear map $\widetilde{A}$ between the corresponding exterior algebras   $\bigwedge \R^{n\ast}$ and $\bigwedge  \R^n$, respectively, by forcing $\widetilde{A}$ to preserve the wedge product. Denote by $A^{(k)}$ the restriction of $\widetilde{A}$ acting  between elements of degree $k$ (this is well--defined, since $\widetilde{A}$ has degree 0). For instance, $A^{(0)}=1$, $A^{(1)}=A$, and $A^{(n)}=\det A$, under obvious identifications. The curious reader may also check that $A^{(n-1)}=A^\sharp$, the cofactor matrix of $A$. On the top of that, each $A^{(k)}$ turns out to be symmetric, that is $A^{(k)}\in S^2 \bigwedge^k \R^n$. This way,  we have  defined the injective map
\begin{eqnarray}
S^2\R^n &\longrightarrow & \p\left( \bigoplus_{k=0}^n  S^2 \bigwedge^k \R^n \right)\, ,\label{eqEQprotoPluck}\\
A &\longmapsto & [(A^{(0)},A^{(1)},\ldots,A^{(n)})]\, ,\nonumber
\end{eqnarray}
and it can be proved that  $\LL_{n,2n}$ is precisely the closure of its image of.\par Interestingly enough, the range of the map \eqref{eqEQprotoPluck} can be made smaller. More precisely, there exists a proper projective subspace $\p V_{\lambda_n}$ which contains the image of \eqref{eqEQprotoPluck} and is minimal with respect to this property. On a deeper level, $\LL_{n,2n}$ should be regarded as a homogeneous manifold of the Lie group $\Sp_{2n}$ and  $V_{\lambda_n}$ as the irreducible $\Sp_{2n}$--representation   realising $\LL_{n,2n}$ as a projective variety in $\p V_{\lambda_n}$ \cite[Section 5.1]{MR3603758}.\par 
We insisted on the fact that $\LL_{n,2n}$ contains the linear space $S^2\R^n$ as an open and dense subset, because this point of view allows to immediately see that
\begin{equation}\label{eqEQprotoTangLG}
T_L\LL_{n,2n}\equiv S^2\R^n\, ,
\end{equation}
for \emph{all} $L\in \LL_{n,2n}$. In other words, the tangent geometry of $\LL_{n,2n}$ is modeled by $n\times n$ symmetric matrices. In fact, the noncanonical identification \eqref{eqEQprotoTangLG} becomes canonical if $\R^n$ is replaced by $L^*$, viz.
\begin{equation}\label{eqEQIsoFund}
T_L\LL_{n,2n}=S^2L^*\, .
\end{equation}
On the   canonical identification \eqref{eqEQIsoFund} alone stands the  rich   geometric approach to $2\Nd$ order PDEs based on contact manifolds. Indeed, since   $\LL_{n,2n}$ is the abstract fibre of $M^{(1)}$, the element $L$ may be thought of as a point of $M^{(1)}$. Hence, in view of the canonical character of the indentification \eqref{eqEQIsoFund}, we as well have
\begin{equation}\label{eqEQIsoFundLOC}
T_LM^{(1)}_p  =S^2L^*\, ,
\end{equation}
for all points $p\in M$. Or, in an equivalent but more abstract way,
\begin{equation}\label{eqEQIsoFundGLOB}
VM^{(1)}=S^2L^*\, .
\end{equation}
Observe that now we have the \emph{vertical bundle} $VM^{(1)}$ and the \emph{tautological bundle} $L$, both over the same base $M^{(1)}$. Formula \eqref{eqEQIsoFundGLOB} explains their global interrelationship, whereby  \eqref{eqEQIsoFundLOC} captures it only on the level of a single fibre. The reader must be extra careful, since the same symbol $L$ denotes both an element of $M^{(1)}$ (as in \eqref{eqEQIsoFundLOC}) and the $n$--dimensional bundle $L\longrightarrow M^{(1)}$, whose fibre at $L$ is, by definition,  $L$ itself (as in \eqref{eqEQIsoFundGLOB}).
\subsection{The symbol of a $2\Nd$ order PDE}\label{subSymb}
If a $2\Nd$ order PDE is understood as a hypersurface \eqref{eqEQ2ndOrdPDE} in $M^{(1)}$, then its vertical bundle $V\E$ is a sub--bundle of $VM^{(1)}$, of codimension 1. Hence, its annihilator $(V\E)^\circ$ is a well--defined one--dimensional sub--bundle (a.k.a. \emph{line bundle}) of $S^2L$, called the \emph{symbol} of the equation $\E=\{F=0\}$. In practice, one can use $F$ to find a (noncanonical) generator of such a line bundle, usually denoted by $S(F)$, viz.
\begin{equation}\label{eqEQSLF}
(V_L\E)^\circ=\Span{S_L(F)}\in \p S^2L\, .
\end{equation}
Equation \eqref{eqEQSLF} intrinsically defines $S(F)$ up to a projective factor. In local coordinates,
the definition of $S_L(F)$ is precisely the one given by  \eqref{eqSIMBOLOprimordiale}, where now the derivatives of $F$ have to be evaluated at the symmetric $n\times n$ matrix corresponding to the Lagrangian $n$--plane $L$ via identification \eqref{eqLagrPlane}. The   reader should not forget that, in spite of the abstract flavour of the definition \eqref{eqEQSLF}  of the symbol, its representative  $S(F)$ is easily computed.\par
Now it is clear what is the $n$--dimensional linear space mentioned in the Introduction (for $n=2$), such that the symbol is a tensor over it. It is precisely $L$. Then   the elements  $\partial_i$ appearing in \eqref{eqSIMBOLOprimordiale} represent a basis of $L$, in compliance with   \eqref{eqLagrPlane}  (indeed, given $L$, it is possible to choose contact coordinates in such a way that $u_{ij}=0$).
%
\subsection{Characteristics of a $2\Nd$ order PDE}\label{subChar}
If solutions of \eqref{eqEQ2ndOrdPDE} are understood as Lagrangian submanifolds $U\subset M$  such that $U^{(1)}\subset M$, then an \emph{initial condition} must be understood as an $(n-1)$--dimensional\footnote{See \cite{MorenoCauchy} for a thorough discussion on the geometry of Cauchy data for nonlinear PDEs.} integral submanifold $\Sigma\subset M$ of the contact EDS $(M,\theta)$. An initial value problem looks now as simple as a pair $(\E,\Sigma)$. A candidate solution $U^1_f$ is then a solution to the initial value problem $(\E,\Sigma)$ if and only if $U^2_f\subset \E$ and also 
\begin{equation}\label{eqEQSigmaContUEffe}
\Sigma\subset U^1_f\, .
\end{equation}
Observe that, once again, the formalism allowed to express everything in terms of set--theoretical inclusions. In local coordinates, $\Sigma\subset U^1_f$ means that the values of $f$ and its $1\St$ derivatives are assigned over an hypersurface in the space $(x_1,\ldots, x_n)$ of the independent variables, thus recovering the familiar picture of an initial condition.\par 
Let $p\in \Sigma$ and $L\in\E$, such that 
\begin{equation}\label{eqEQLcontieneTSigma}
 T_p\Sigma\subset L \, .
\end{equation}
 Observe that $L$ is always of the form $L=T_pU^1_f$ for some candidate solution $U^1_f$, so that \eqref{eqEQLcontieneTSigma} is exaclty the infinitesimal version of \eqref{eqEQSigmaContUEffe} at $p$. Accordingly, we should think of $L$ as an \emph{infinitesimal solution} of $\E$, and we should think of $ T_p\Sigma$ as an \emph{infinitesimal initial datum} (at $p$).  The  above--introduced symbol $S(F)$ allows to answer the following question: is $L$ the unique infinitesimal solution passing through that infinitesimal initial datum? The answer is simple and operative. Regard $T_p\Sigma$ as a line in $L^*$, i.e., as a point of $\p L^*$. Since    $\p L^*$ contains the quadric hypersurface of equation $S_L(F)=0$ (see \eqref{eqEQSLF}), there are only two options: either   $T_p\Sigma$ belong to that quadric, or it doesn't. In the first case,   the answer is negative and $T_p\Sigma$ is a characteristic hyperplane.\par 
 Due to its importance, the quadric $S_L(F)=0$ is called the \emph{characteritic variety} of $\E$ at $L$. A line $\ell\in\p L^*$ (that is, an hyperplane in $L$) is \emph{characteristic} if it belongs to the characteristic variety. In conclusion, we have obtained a natural geometric picture of ill--defined (infinitesimal) initial value problems: the characteristic lines are the normals to the tagent hyperplanes to those initial data for which the Cauchy--Kowalevskaya theorem fails in uniqueness \cite[Section 1.2]{Vitagliano2013}.  \par 
It is worth observing that any line $\ell\in\p L_0^*$ (i.e., an hyperplane $\ell^\circ\subset L_0$) determines a \emph{curve}
\begin{equation}\label{eqEQcurvEllUno}
\ell^{(1)}:=\{L\in M^{(1)}\mid L\supset\ell^\circ\}
\end{equation}
in $M^{(1)}$ passing through $L_0$. Then the reader may verify that $\ell$ is characteristic at $L_0$ (i.e., it lies in the quadric $S_{L_0}(F)=0$) if and only if the curve $\ell^{(1)}$ is tangent to $\E$ at $L_0$. This alternative interpretation immediately  leads us to       another   definition. If the curve  $\ell^{(1)}$ is entirely contained into $\E$, then $\ell$ is a \emph{strong characteristic} line \cite[Section 3.1]{MR2985508}.
\subsection{The test $S^2(F)\approx S(F)$}
We are now ready to motivate the rank--4 tensor $S^2(F)$ introduced by \eqref{eqEQS2Fprimordiale}, and to explain the condition \eqref{eqCEprimordiale}.\par
The quickest way to get to $S^2(F)$ inevitably sacrifices     coordinate--independence. The construction goes as follows. Regard the components of the tensor \eqref{eqSIMBOLOprimordiale} as functions on $M^{(1)}$, and replace each of them by its own symbol. The result
\begin{equation}
S\left( \frac{\partial  F}{\partial u_{ij} }\right)\partial_i\odot\partial_j\, ,
\end{equation}
interpreted as a homogeneous polynomial of degree 4, is precisely \eqref{eqEQS2Fprimordiale}. Now, both $S(F)$ and $S^2(F)$ are elements of the polynomial algebra $S^\bullet L$, and as such it is legitimate to ask whether the latter is proportional to the former. That is, condition \eqref{eqCEprimordiale} makes sense. So, we can use it to single out a nontrivial class of $2\Nd$ order PDEs.
\begin{definition}\label{defBucchinariellaCheMiLevaDaiPasticci}
The equation $\E$ \emph{passes the test $S^2(F)\approx S(F)$} if and only if the condition \eqref{eqCEprimordiale} is satisfied by some $F$  such that  $\{F=0\}=\E$.
\end{definition}
It is worth stressing that it is the \emph{equation} $\E$ cut out by $F$ which passes the test $S^2(F)\approx S(F)$, and not $F$ itself (see \cite[Proposition 3.6]{MR3603758}).  The collection of all $2\Nd$ order PDEs passing this test is precisely  the nontrivial class of $2\Nd$ order PDEs introduced by G. Boillat \cite{MR1139843}, following P. Lax's original intuition that there should exists a class of nonlinear PDEs ``not genuinely nonlinear''.   By construction, this class is automatically  invariant under the symmetry group of the theory, that is the group of all contactomorphisms of $M$.  This is the class of completely exceptional (scalar, nonlinear) $2\Nd$ order PDEs.

\section{Multidimensional Monge--Amp\`ere equations}\label{sec2}
Since a linear PDE manifestly passes the test $S^2(F)\approx S(F)$,   it seems natural to suspect that the class of  completely exceptional PDEs is the closure, under the action of the group of contactomorphisms of $M$, of the class of linear PDEs---we call such a closure the class of \emph{linearisable} PDEs \cite{MR1606791}. But life is slightly harder (and interesting) than that: linearisable PDEs form  a proper subclass in the class of completely exceptional PDEs. The purpose of this section is to show that the latter coincides with the class of (multidimensional) Monge--Amp\`ere equations on $M$.\par 
A Monge--Amp\`ere equation on $M$ can be seen as an EDS on $M$, namely
\begin{equation}\label{eqEQMAEDS}
(M,\{\theta,\omega\})\, ,
\end{equation}
where $\omega\in\Omega^n(M)$ is an $n$--form on $M$. Recall that the set of ($n$--dimensional) integral elements of the contact EDS $(M,\theta)$ is the bundle $M^{(1)}\longrightarrow M$. Hence, the set of the integral elements of the EDS \eqref{eqEQMAEDS} will be a proper subset of $M^{(1)}$. In fact, it will be a hypersurface, henceforth denoted by $\E_\omega$, that is a $2\Nd$ order PDE according to our understanding (see Section \ref{subPDEs}). Since   we are looking for $n$--dimensional  integral elements, the condition $d\omega|_L\equiv 0$ is vacuous, that is we are just imposing the unique condition $\omega|_L\equiv 0$ on the integral elements of $(M,\theta)$. This explains the codimension one. The  definition of Monge--Amp\`ere equations as hypersurfaces of the form $\E_\omega$ was given in 1978 by  V. Lychagin \cite{Kushner2009}. \par 
Alternatively, Monge--Amp\`ere equations can be defined as \emph{hyperplane sections} of the fibres of $M^{(1)}$. Recall from Section \ref{subsecLG} that each fibre $M^{(1)}_p$ of the bundle $M^{(1)}\longrightarrow M$ is a projective variety in $\p V_{\lambda_n}$. Then one can define Monge--Amp\`ere equations as the codimension--one sub--bundles $\E\subset M^{(1)}$, such that each fibre $\E_p$ is the intersection of $M^{(1)}_p$ with an hyperplane of $\p V_{\lambda_n}$. Such intersection is called an hyperplane section of the Lagrangian Grassmannian $M^{(1)}_p$ .\par 
This definition is by no means new, it is only formulated in an abstract way. Bearing in mind the embedding \eqref{eqEQprotoPluck}, we easily see that a hyperplane section of $M^{(1)}_p$ is locally given by a (unique) \emph{linear relation} between the minors of the Hessian matrix $A:=\|u_{ij}(p)\|$ of $u$ computed at $p$ (a simple exercise of multilinear algebra).  By letting the point $p$ vary on $M$, we obtain the coordinate expression of the (left--hand side of) Monge--Amp\`ere equations, viz.
\begin{equation}\label{eqEQcoordMA}
F=B_0+B_1^{ij}u_{ij}+B_2^{\ldots}(2\times 2\textrm{ minors})+\cdots+B^{ij}_{n-1}u^\sharp_{ij}+B_n\det\| u_{ij}\|\, ,
\end{equation}
where now $B_0,B_1^{ij},\ldots, B_n$ are functions on $M$. The local expression \eqref{eqEQcoordMA} is the classical way multidimensional Monge--Amp\`ere equations are introduced. Its equivalent  interpretation in terms of hyperplane sections came later (see, e.g., \cite[Section 3.2]{MR2985508}).  \par 
Formula \eqref{eqEQcoordMA} shows that linear equations are in particular Monge--Amp\`ere equations (enough to set $B_2=B_3=\cdots=B_n=0$). Less evident is that, if $F$ is of the form \eqref{eqEQcoordMA}, then $S^2(F)\approx S(F)$.  
The final result of this section tells precisely what are the solutions to the equation \eqref{eqCEprimordiale}. This result was firstly  proved by G. Boillat \cite{MR1139843} and re--examined  recently    \cite{MR3603758}.  
\begin{theorem}
The   $2\Nd$ order PDE $\E=\{F=0\}$ passes the test $S^2(F)\approx S(F)$ (see  Definition \ref{defBucchinariellaCheMiLevaDaiPasticci}) if and only if $\E$ is a Monge--Amp\`ere equation.
\end{theorem}
From an abstract standpoint---that is by giving up any interpretation of the objects at play in terms of (nonlinear) PDEs---the equation $S^2(F)\approx S(F)$ is the ``flatness condition'' for a (family of) hypersurface in the Lagrangian Grassmannian. By reprising the first (very trivial) equation of this paper \eqref{eqEqBanLin}, we may observe that 
\begin{equation}\label{eqEQTrivial}
\frac{\partial^2 q}{\partial z_a\partial z_b}=0
\end{equation}
is the condition for the hypersurface $\{q=0\}$ in $\p V_{\lambda_n}$ to be a hyperplane (the $z_a$'s are projective coordinates and $q$ is a homogeneous polynomial). Since $\LL_{n,2n}\subset \p V_{\lambda_n}$, it is natural to ask how the (very trivial) condition \eqref{eqEQTrivial} looks like if one knows only the restriction $F:=q|_{\LL_{n,2n}}$. The (nontrivial\footnote{See \cite[Section 5]{MR3603758} for a taste of the techniques involved.}) answer is: precisely $S^2(F)\approx S(F)$. The subtle point is that equation \eqref{eqEQTrivial} recognises when $q$ is linear in the coordinates $z_a$'s, whereas $S^2(F)\approx S(F)$ recognises when $F$ is \emph{linear in the minors} of the symmetric $n\times n$ matrix $\|u_{ij}\|$, that is, when $F$ is of the form \eqref{eqEQcoordMA}. Observe that, as a function of the \emph{entries} of $\|u_{ij}\|$, a function $F$ as in \eqref{eqEQcoordMA} is (nonhomogeneously) of degree $\leq n$.\par 
Now we start thinking at a \emph{solution} $F$ of the equation $S^2(F)\approx S(F)$ as the left--hand side of a $2\Nd$ order PDE itself. We ask ourselves what makes $F$ so special amongst all possible left--hand sides of $2\Nd$ order PDEs. The answer to this question is in fact the original 1954 observation by P. Lax: there exist certain  (scalar, generally nonlinear) $2\Nd$ order PDEs that are characterised by an ``exceptional  behaviour'' of their solutions. In continuity with P. Lax's work, G. Boillat later called them ``completely exceptional'', and now we know that this class coincides precisely with the class of Monge--Amp\`ere equations \cite{MR3603758}.\par 
The modern geometric reinterpretation of P. Lax's class of equations, succinctly captured by the condition $S^2(F)\approx S(F)$, by no means diminishes the value of his work. On the contrary it shows that inside this class, whose existence follows from abstract representation--theoretic arguments, one finds several different PDEs, an yet all these PDEs share an important and physically meaningful property: their solutions display a sort of ``linear behaviour''. This property, and its link with the condition  $S^2(F)\approx S(F)$, will be explained in the last section below. 
\section{Discontinuity waves, shock waves and completely exceptional PDEs}\label{sec3}
From now on we assume $n=2$ and we deal only with hyperbolic PDEs. This is necessary to shorten the distances between the original P. Lax/G. Boillat's idea of a completely exceptional $2\Nd$ order PDE and the modern interpretation of Monge--Amp\`ere equations. 
\subsection{Hyperbolic $2\Nd$ order PDEs all whose characteristics are strong}\label{subHypPDE}
Let then $\E=\{F=0\}$ be a   2--dimensional hyperbolic $2\Nd$ order PDE. In this case, $M$ is a 5--dimensional contact manifold, and the generic fibre of $M^{(1)}$ is the 3--dimensional Lagrangian Grassmannian $\LL_{2,4}$.  The irreducible $\Sp_4$--module $V_{\lambda_2}$ is 5--dimensional, and hence $\p V_{\lambda_2}=\p^4$. It is known that $\LL_{2,4}$ is a smooth quadric hypersurface (called Lie quadric) and a conformal manifold as well \cite{MR2876965}. \par
A remarkable (and general) property of Lagrangian Grassmannians is that the curves $\ell^{(1)}$ (cf. \eqref{eqEQcurvEllUno}) are straight lines in the projective space $\p V_{\lambda_n}$.   Suppose that $\ell\in \p L^*$ is a characteristic of $\E$ at $L\in\E$. Then, by definition (see the end of Section \ref{subChar}) the curve  $\ell^{(1)}$ is tangent to $\E$ at $L$.  Now we use the hyperbolicity of $\E$, that is, the fact that $S_L(F)$ has two distinct roots. The characteristic variety (see Section \ref{subChar}) is then a degenerate quadric, namely the union $\ell_1\cup\ell_2\subset L^*$ of two lines. Since the polynomial $S_L(F)$ depends on $L\in\E$, the dependency of $\ell_i$  on $L$ should be stressed, and from now on we write $\ell_{i,L}$, and let $i$ be either 1 or 2. Consequently, to any infinitesimal solution $L\in\E$, we can can associate a line $\ell_{i,L}^{(1)}$ in $\p^4$, i.e., we have a sort of Gauss map defined on $\E$ and taking its values in the Grassmannian of all projectives lines in $\p^4$. As such, the derivative of the map $L\longmapsto \ell_{i,L}^{(1)}$ can be computed---and this is the crucial point---\emph{along the curve  $\ell_{i,L}^{(1)}$ itself at the point $L$}. Without loss of generality, let us assume that $\ell_{i,L}^{(1)}$ passes through $L$ for $t=0$, and let $\dot{\ell}_{i,L}^{(1)}(0)$ denote its velocity thereby.
\begin{theorem}\label{thBellino}
 The equation
 \begin{equation}\label{eqEQpseudogeodetica}
\left(  \dot{\ell}_{i,L}^{(1)}(0)\right)(\ell_{i}^{(1)})=0\, 
\end{equation}
is satisfied for all $L\in\E=\{F=0\}$ if and only if $S^2(F)\approx S(F)$.
\end{theorem}
Theorem \ref{thBellino} is the last step towards    our main purpose, which  is to prove that  the test $S^2(F)\approx S(F)$ is equivalent to original Lax's condition of complete exceptionality. Indeed, as we will show in the next subsection, in appropriate local coordinates the condition \eqref{eqEQpseudogeodetica} takes precisely the form of Lax's condition.\par
The proof of Theorem \ref{thBellino} is not hard and can be found in \cite[Corollary 3.8]{MR3603758}. Here we can   provide  an intuitive explanation. The lines $\ell_i$ are just the roots of the polynomial $S(F)$, but in our hypothesis of hyperbolicity knowing the polynomial is exactly the same as knowing both its roots (up to a projective factor, but this does not matter, cf. \eqref{eqEQSLF}). The condition \eqref{eqEQpseudogeodetica} means that these roots do not vary along some special directions determined by the roots themselves. In a sense, the condition $S^2(F)\approx S(F)$ is dual to \eqref{eqEQpseudogeodetica}: instead of taking a particular derivative of the roots, we took the whole differential of the corresponding polynomial, and then we imposed its vanishing on the zero locus of the polynomial itself, i.e., the same special directions as before.  It is worth stressing how the reformulation of \eqref{eqEQpseudogeodetica} as $S^2(F)\approx S(F)$ freed us from the necessity of a complete decomposable polynomial. Hyperbolicity condition is however  indispensable for interpreting $S^2(F)\approx S(F)$  in terms of behaviour of solutions.\par
Before recalling P. Lax's original observation, let us comment further on the condition \eqref{eqEQpseudogeodetica}, stressing a detail which could not be appreciated when  reformulated as $S^2(F)\approx S(F)$.   The 
equation  \eqref{eqEQpseudogeodetica} was deliberately laid down in a form which is reminiscent of the equation of geodesics---curves whose speed is constant along the curves themselves. Indeed, if \eqref{eqEQpseudogeodetica} vanish identically, it means that the family  $L\longmapsto \ell_{i,L}^{(1)}$ of   projective lines in $\p^4$ passes with speed zero through the point $L_0$ of $\ell_{i,L_0}^{(1)}$,  which is now understood   as a (vertical) curve tangent to $\E$. In other words, for a fixed $L_0\in \E$, the curve  $\ell_{i,L_0}^{(1)}$  coincides, for all points  $L\in \ell_{i,L_0}^{(1)}$, with the curve $\ell_{i,L}^{(1)}$  associated to $L$. As such, $\ell_{i,L_0}^{(1)}$ is tangent  to $\E$ in all its points and hence is \emph{entirely contained} into $\E$. We have therefore proved that, if \eqref{eqEQpseudogeodetica} is identically satisfied for all $L\in\E$, then \emph{any characteristic of $\E$ is a strong characteristic}. The converse is also true, and even easier to prove \cite[Theorem 5.9]{MR2985508}.\par
The conclusion of this subsection is that, for hyperbolic PDEs, the condition  $S^2(F)\approx S(F)$ can be recast in a somewhat more tangible form: a PDE satisfies it if and only if all its characteristics are strong characteristics. We see then how an abstractly defined class of PDEs is characterised by a special behaviour of their characteristics, which are ultimately linked to questions of existence and uniqueness of solutions. It is only in P. Lax's original formulation that one can see how actual solutions behave.
\subsection{The ``not genuinely nonlinear'' nonlinear PDEs in the sense of P. Lax}
For $n=2$, the symbol \eqref{eqSIMBOLOprimordiale}, computed at $L\in\E$, reads
\begin{equation}
S(F)=F_{u_{xx}}\partial_x^2+F_{u_{xy}}\partial_x\partial_y+F_{u_{yy}}\partial_y^2\, .
\end{equation}
Recall that $\partial_x$ and $\partial_y$ are the generators of  $L$, and that  $S(F)$ is a quadratic polynomial on $L^*$ (see Section \ref{subSymb}).   A line\footnote{As the reader may have guessed, $dx$  and $dy$ are dual to $\partial_x$ and $\partial_y$, respectively.} $\ell=\Span{\xi dx+\eta dy}\in\p L^*$ is characteristic iff
\begin{equation}
F_{u_{xx}}\xi^2+F_{u_{xy}}\xi\eta+F_{u_{yy}}\eta^2=0\, ,
\end{equation}
that is (assuming $\xi\neq 0$),
\begin{equation}
F_{u_{xx}} +F_{u_{xy}}\lambda+F_{u_{yy}}\lambda^2=0\, ,
\end{equation}
where $\lambda=\eta/\xi$. The function\footnote{We often neglect stressing the dependency upon  $L\in\E$ of the quantities at play. } $\lambda$  is the \emph{characteristic speed}, according to P. Lax. By the hyperbolicity assumption, there are two functions $\lambda_i$, $i=1,2$, such that
\begin{equation}\label{eqEQElleIElle}
\ell_{i,L}=\Span{dx+\lambda_i(L)dy}\, ,\quad i=1,2\, .
\end{equation}
Recall (see Section \ref{subHypPDE} above) that the union   $\ell_{1,L}\cup\ell_{2,L}$ is precisely the zero set of $S_L(F)$. 
Consider then the curves $\ell_{i,L}^{(1)}$ passing through $L$ at time 0   and compute their speed at 0. 
It is a simple exercise \cite[Proposition 2.5]{MR3603758} to prove that
\begin{equation}\label{eqEQDotElleIElle}
\dot{\ell}^{(1)}_{i,L}(0)=\partial_{u_{xx}}+\lambda_i(L)\partial_{u_{xy}}+\lambda^2_i(L)\partial_{u_{yy}}\, .
\end{equation}
Formulae \eqref{eqEQElleIElle} and \eqref{eqEQDotElleIElle} allows us to write down the condition \eqref{eqEQpseudogeodetica} in terms of the characteristic speed. Indeed, locally, $\ell_i^{(1)}$ is completely characterised by the characteristic speed $\lambda_i$ and hence, instead of asking that the derivative of the former be zero, we may require the derivative of the latter be zero. In other words, \eqref{eqEQpseudogeodetica} is equivalent to 
\begin{equation}
\left(  \dot{\ell}_{i,L}^{(1)}(0)\right)(\lambda_{i} )=0\, ,
\end{equation}
that is
\begin{equation}\label{eqEQ14}
(\lambda_i)_{u_{xx}}+\lambda_i(\lambda_i)_{u_{xy}}+\lambda^2_i(\lambda_i)_{u_{yy}}=0\, ,\quad i=1,2\, .
\end{equation}
Equation \eqref{eqEQ14} is the local counterpart of $S^2(F)\approx S(F)$, as can be found in \cite[Equation (14)]{MR3603758} or \cite[Equation (5)]{MR1292999}.  Its derivation, as presented here, is however antihistorical, since \eqref{eqEQ14} was    formulated  decades   before $S^2(F)\approx S(F)$. Observe the local and coordinate--dependent nature of \eqref{eqEQ14}, as  opposed to the intrinsic character of $S^2(F)\approx S(F)$. However, in the formulation \eqref{eqEQ14}  there enter   the characteristic speeds $\lambda_i$, and these allow for a tangible interpretation of complete exceptionality in terms of behaviour of solutions. \par
In his 1954 paper \cite{MR0068093} P. Lax was interested in the behaviour of solutions of  certain ($1\St$ order hyperbolic, systems of) PDEs. Though he suspected and postulate the existence of a distinguished class of such PDEs, the problem whether or not a special ``linear behaviour'' of the solutions to a nonlinear PDE could  unambiguously define a class of  PDEs was marginal to him.  He was mainly interested in  the phenomenon of development of shocks out of ``weak discontinuities'', and his key remark was that discontinuities propagates along characteristics (cf. Theorem \ref{thBellino}). \par
\begin{figure}
 \epsfig{file=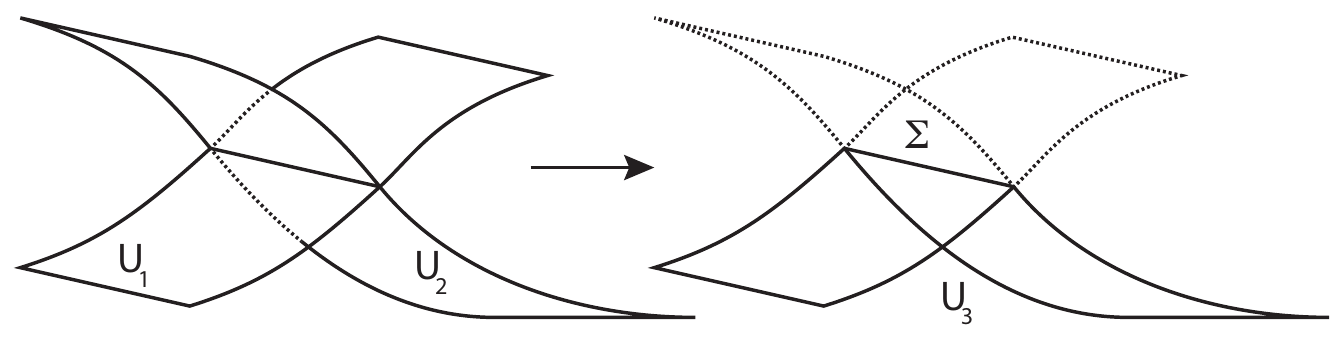,width=\textwidth}
 \caption{If two smooth solutions, say $U_1$ and $U_2$ are combined, a nondifferentiable  solution $U_3$ is obtained (a weak discontinuity). By construction, the locus $\Sigma$ where $U_3$ is singular corresponds to     characteristic initial data.\label{figDiscWaw}}
\end{figure}
\begin{figure}
 \epsfig{file=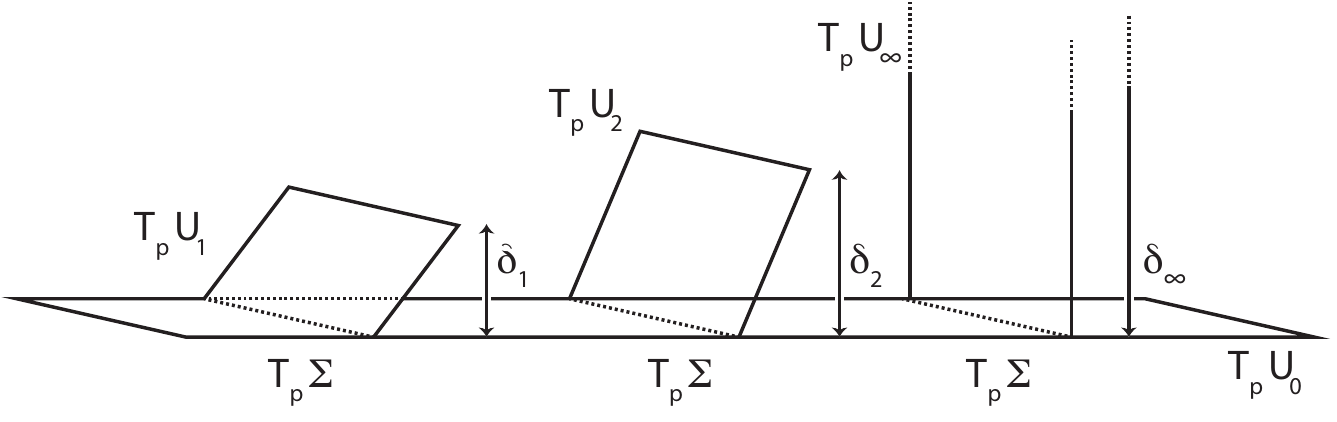,width=\textwidth}
 \caption{If  the ``steady state solution'' $U_0$ is glued with another solution $U_i$ along $\Sigma$, then   a weak discontinuity   is obtained. The ``jump'' $\delta$ measure that the discrepancy between the corresponding tangent spaces. If $\delta$ is free to grow, then a weak discontinuity   may  continuously evolve into a shock---a typical nonlinear phenomenon.\label{FigJump}}
\end{figure}

In the present  present context, his original remarks may be recast as follows. Let us interpreted  a small difference $ {\ell}^{(1)}_{i,L}(\varepsilon)-{\ell}^{(1)}_{i,L}(0)$ as the \emph{jump} in the $1\St$ derivatives of a continuous and almost everywhere smooth solution of $\E$---what he called a \emph{weak discontinuity}. Let us call this difference simply $\delta$, as in \cite{MR1292999}.   According to the understanding of a characteristic as an (infinitesimal) ill--defined initial value problem, a weak discontinuity  can be obtained by gluing two pieces belonging to different solutions passing through the same (characteristic) initial datum (see Figure \ref{figDiscWaw}), but with different values of ``normal derivatives''. The discrepancy between the tangent spaces $L={\ell}^{(1)}_{i,L}(0)$ and $L_\varepsilon:= {\ell}^{(1)}_{i,L}(\varepsilon)$ of the two solutions passing through the same initial datum $\ell_{i,L}^\circ$ is measured (infinitesimally) by $\delta$---the  ``jump'', see Figure \ref{FigJump}. The ``infinitesimal version'' of $\delta$ is precisely the tangent vector $\dot{\ell}^{(1)}_{i,L}(0)$ which, in view of \eqref{eqEQDotElleIElle}, is fully described by the value $\lambda_{i}(L)$.\par

\begin{figure}
 \epsfig{file=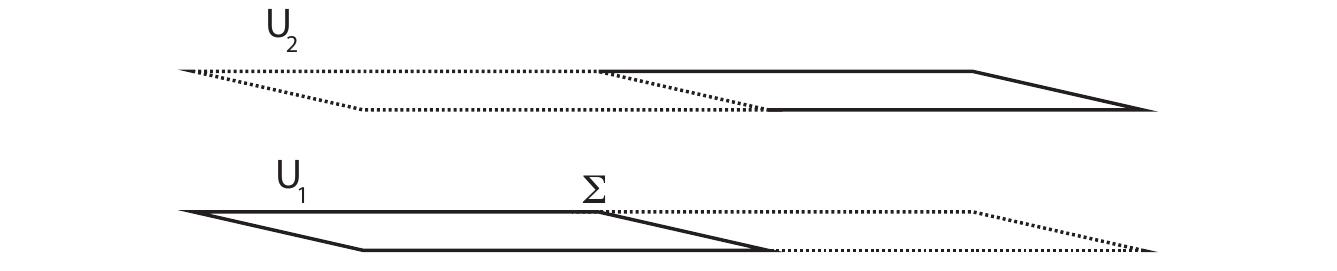,width=\textwidth}
 \caption{A shock is a discontinuous solution where the very value of the solution experiences a jump along $\Sigma$.\label{FigShock}}
\end{figure}
Now we can rephrase the definition of completely exceptional PDEs (implicitly) present in P. Lax's paper \cite{MR0068093}: \emph{a completely exceptional PDEs is a (nonlinear) PDEs whose weak discontinuities   never evolve into shocks in a finite time}.\footnote{By \emph{time} here we mean a field transversal to the wavefronts, that is characteristic surfaces.} A \emph{shocks}  is a solution almost everywhere smooth but not even continuous along characteristics (see Figure \ref{FigShock}). Intuitively, this amounts at requiring the ``jump'' $\delta$ to be constant along the characteristics themselves, which immediately translates into \eqref{eqEQ14}, which in turn is equivalent to \eqref{eqEQpseudogeodetica}, which is a particular case of the condition $S^2(F)\approx S(F)$ for hyperbolic PDEs. So, all definitions are equivalent on their common ground. 

\subsection{Concluding remarks}
We have proved that the class of $2\Nd$ order (scalar, nonlinear) PDEs passing the test $S^2(F)\approx S(F)$ is an enlargement of the class of (quasi)linear PDEs. In oder to achieve this enlargement is it however not enough to apply all possible contactomorphism to the (quasi)linear PDEs, because this generates a proper sub--class. The correct interpretation of the PDEs $\E=\{F=0\}$ such that $S^2(F)\approx S(F)$ is as those ``nonlinear PDEs displaying a linear behaviour in their solutions'', meaning that weak discontinuities never evolve into shocks---a typical feature of linear PDEs. It is truly remarkable that P. Lax's class of completely exceptional PDEs turned out to coincide with the class of Monge--Amp\`ere equations, which are defined in purely algebraic terms as (families of) hyperplane sections of the Lagrangian Grassmannian.

\section*{Acknowledgements}
This note   is an adaptation of the poster titled ``A representation--the\-o\-ret\-ic characterisation of completely exceptional second--order PDEs'', which was based on the joint paper \cite{MR3603758} with J. Gutt and G. Manno  and was presented by the author at the $50\Th$ Seminar ``Sophus Lie''  in Bedlewo, Poland, 25 September -- 1 October 2016.   The  research of the author has been partially supported by the Marie Sk\l odowska--Curie fellowship SEP--210182301 ``GEOGRAL", and has also been partially    founded by the  Polish National Science Centre grant
under the contract number 2016/22/M/ST1/00542.  The author thanks the organisers of the Seminar ``Sophus Lie'' for their excellent job and the anonymous referee for pointing out a clumsy oversight in the preliminary version of this paper. The author  is a  member of G.N.S.A.G.A. of I.N.d.A.M. 


\def\polhk#1{\setbox0=\hbox{#1}{\ooalign{\hidewidth
  \lower1.5ex\hbox{`}\hidewidth\crcr\unhbox0}}}

\end{document}